\documentclass{article}
\usepackage{amsmath, amssymb, mathrsfs, eucal, latexsym}
\usepackage{color}

\newtheorem{theorem}{Theorem}[section]
\newtheorem{prop}[theorem]{Proposition}
\newtheorem{lemma}[theorem]{Lemma}
\newtheorem{remark}[theorem]{Remark}

\newtheorem{example}[theorem]{Example}
\newtheorem{corollary}[theorem]{Corollary}

\def\bR{\mathbb{R}}
\def\bZ{\mathbb{Z}}
\def\bE{\mathbb{E}}
\def\bH{\mathbb{H}}
\def\bQ{\mathbb{Q}}
\def\bN{\mathbb{N}}
\def\bD{\mathbb{D}}
\def\bL{\mathbb{L}}

\def\cA{\mathcal{A}}
\def\cB{\mathcal{B}}
\def\cC{\mathcal{C}}
\def\cD{\mathcal{D}}
\def\cE{\mathcal{E}}
\def\cF{\mathcal{F}}
\def\cG{\mathcal{G}}
\def\cH{\mathcal{H}}
\def\cI{\mathcal{I}}
\def\cK{\mathcal{K}}
\def\cP{\mathcal{P}}
\def\cM{\mathcal{M}}
\def\cO{\mathcal{O}}
\def\cR{\mathcal{R}}
\def\cS{\mathcal{S}}
\def\cT{\mathcal{T}}
\def\cL{\mathcal{L}}
\def\cU{\mathcal{U}}

\begin{document}
	
	\title{Strong consistency and optimality for generalized estimating equations with stochastic covariates \\
	}
	\date{\today}
	\author{Laura Dumitrescu \textsuperscript{a} \footnote{\noindent Corresponding author. E-mail: laura.dumitrescu@vuw.ac.nz }
		\and Ioana Schiopu-Kratina \textsuperscript{b}
	}
	\maketitle
	\noindent \textsuperscript{a}{Victoria University of Wellington, New Zealand} \\
	\noindent \textsuperscript{b}{Department of Mathematics and Statistics, University of Ottawa, Canada}
	\begin{abstract}
		\noindent In this article we study the existence and strong consistency of GEE estimators, when the generalized estimating functions are martingales with random coefficients. Furthermore, we characterize estimating functions which are asymptotically optimal.
		\vspace{1mm}
		
	\end{abstract}
	
	
	
	\section{Motivation} 
	There has been significant interest and work dedicated to constructing theoretical and methodological tools for the analysis of longitudinal data. A popular estimation method, based on generalized estimating equations (GEEs), was introduced by Liang and Zeger (1986). Several asymptotic results have appeared since this seminal paper. The literature on this topic handles deterministic regressors, while, to the best of our knowledge, asymptotic results for GEEs with stochastic regressors have not been discussed. In related topics, results on models with stochastic regressors were considered in Lai and Wei (1982) and Chen et al. (1999), and strong convergence of estimators was established for the linear and generalized linear models, respectively. The property of strong consistency of estimators is often of interest, for example, in stochastic control problems, where the estimate is sequentially updated at each stage. Furthermore, applications of regression models with stochastic regressors include empirical modelling of economic learning models, spatial statistics or neural network models (see Lai, 2003 for a list of references). 
	
	In the context of GEEs for longitudinal data, weak consistency holds even under misspecification of the correlation structure (see Liang and Zeger, 1986 and Xie and Yang, 2003). Conditions for asymptotic optimality were given in Jiang et al. (2007) and Balan et al. (2010).  However, the asymptotic theory cannot be directly translated to the case of stochastic regressors, and our aim here is to formulate asymptotic results which are appropriate for longitudinal models with non-deterministic regressors.

	\section{Our contribution}
	We extend results on the existence, strong consistency and asymptotic optimality of generalized estimation equations (GEEs) estimators to accommodate stochastic regressors in the GEEs. This extension has important consequences, which allow us to give our double-indexed response data two major interpretations, and thus enlarge their domain of applicability. 
	
	Let ${\bf y}_i^T=(y_{i1}, \ldots, y_{im_i})$ represent a response vector of dimension $m_i \le m$, for each $i \ge 1$. The index $i$ may identify an individual, with $m_i$ the number of occasions the response variable is recorded on this individual. Each $i$ is then a cluster of correlated observations, and we may assume that the clusters (the individuals) are independent, $i \ge 1$. This framework is appropriate for longitudinal surveys and constitutes our first interpretation of data. The second interpretation, which we describe next, is only possible when the regressors are stochastic. We view the response vector as a cluster of data evolving in time, which is indexed by $i \ge 1.$ As a function of time, we are dealing with a sequence of martingale transforms, which model a form of dependence of data over time. The vector of data at time $i,$ ${\bf y}_i,$ may represent a family of individuals which are part of a specific neighbourhood. Under this set-up, we may also consider a vector-valued time series with correlated components. For instance, we may be interested in the evolution in time of a portfolio, with options that do not respond independently to the state of the market. In all these instances, there exists a within-cluster positive correlation, which need not be specified.
	
	Our assumption that the cluster sizes are bounded is not unreasonable. In practice, the number of observations recorded on each individual is limited. Likewise, when geographic clusters become too large, they are usually split to decrease the interviewers' load, and one can easily identify an ``old'', surveyed cluster with a newly defined one.  While in Balan et al. (2010) the cluster sizes were all equal to $m$, here they may vary, which is a more realistic assumption, since we often deal with missing observations. 
	
	A GEE estimator is defined as an implicit solution of a generalized estimating equation, and one seeks a sequence of such estimators, which convergences to the true value of the regression parameter and has a minimal asymptotic variance. The GEEs considered here are quite general, in that both the covariate matrices ${\bf X}_i$ and the ``working correlation matrices'' ${\bf R}_i^*,$ $i \ge 1$ are stochastic. The latter matrices may also depend on the main regression parameter.  
	
	
	The first result that we present in Section \ref{optim} is a characterization of asymptotic quasi-score sequences of functions within the class $\mathcal{H}_n$ of estimating functions introduced in Balan et al. (2010), which are square integrable martingale transforms. The asymptotic quasi-score functions are asymptotically optimal GEEs (defined in Heyde, 1997) in the sense that the asymptotic confidence interval defined by the sequence of the GEE estimators associated with these optimal GEEs is minimal. While the class $\mathcal{H}_n$ is large enough to accommodate stochastic regressors, in Balan et al. (2010) we did confine ourselves to characterizing optimal estimating functions with constant regressors, a situation that is more appropriate to the study of longitudinal data with independent clusters. In Section \ref{optim} we show that the class of asymptotic quasi-scores in $\mathcal{H}_n$ is quite large. Not only it includes estimating functions with random coefficients, it also allows for slight misspecification of the regressors when the functions that define the marginal models in the GEEs are continuous (see Corollary \ref{corl}). We expand on this particular type of sensitivity analysis in Section \ref{optim}. 
	
	Our main result on strong convergence is Theorem \ref{strong-consistency-th}, which generalizes Theorem 4.2 of Balan et al. (2010), and thus Theorem 7 of Xie and Yang (2003) to the case of stochastic regressors. Our findings are akin to those published in Lai and Wei (1982) on the topic of strong convergence of the least square estimators. The applications to dynamic systems as discussed in Section 3 of Lai and Wei, (1982) are also covered by our GEE set-up. Essentially, Theorem \ref{strong-consistency-th} is a fixed point theorem, and the use of SLLNs for martingales with a proper normalization is quite important. Our choice of the normalizing factor and the hypotheses of Theorem \ref{strong-consistency-th} ensure that the convergence in the SLLNs is faster than the convergence in the equicontinuity of the normalized derivative. The adequate normalizing random variable is the largest eigenvalue of the conditional variance. The measurability properties of the regressors ensure that our GEEs are transform martingales, and thus martingale strong convergence results apply. The technique we use to prove Theorem \ref{strong-consistency-th} was first introduced in Xie and Yang (2003), which, as remarked in Wang (2011), is appropriate for obtaining strong convergence results. 
	
	Our second result on strong convergence is Theorem \ref{consistency}, in which we break-down the hypotheses of Theorem \ref{strong-consistency-th} into conditions that are easier to verify and facilitate the comparison with previous results. We adapted the technical results in Balan et al. (2010), namely Lemmas 4.7-4.9, to the case of stochastic regressors and random norming, and eliminated or weakened some of the conditions of Theorem 4.14 in Balan et al. (2010).        
	
	\section{Model assumptions}

	Let $\boldsymbol{\beta} \in \mathcal{T} \subset \mathbb{R}^p$ be a possible value of the regression parameter, with $\boldsymbol{\beta}_0 \in \mathcal{T}$ the (unknown) true value of the parameter.  We assume that each ${y}_{ij}$, $j = 1, \ldots, m_i$, $i>1$ is a random variable in the probability space, $(\Omega, \mathcal{F}, P_{\boldsymbol{\beta}})$, with $\boldsymbol{\beta} \in \mathcal{T} \subset \mathbb{R}^p$ , which is equipped with an increasing sequence of sigma-fields, $\mathcal{F}_i$, and that the vector ${\bf y}_i^T$  is $\mathcal{F}_i$-measurable, $i \ge 1$. Let ${\bf X}_i=({\bf x}_{i1}, \ldots, {\bf x}_{im_i})^T$, ${\bf x}_{ij} \in \mathbb{R}^p$, $j=1, \ldots, m_i$ be stochastic regressors (explanatory variables), where each ${\bf x}_{ij}$ is $\mathcal{F}_{i-1}$-measurable, $i \ge 1$. Our model assumptions are 
	\begin{eqnarray*}
		\textrm{E}_{\boldsymbol{\beta}}(y_{ij}| \mathcal{F}_{i-1})  &=& \mu(\textbf{x}_{ij}^T\boldsymbol{\beta}) = \mu_{ij}(\boldsymbol{\beta}),  \\
		\textrm{Var}_{\boldsymbol{\beta}}(y_{ij}| \mathcal{F}_{i-1})&=& \mu'(\textbf{x}_{ij}^T\boldsymbol{\beta}) = \sigma_{ij}^2(\boldsymbol{\beta}),
	\end{eqnarray*}
	\noindent where $\mu$ is a link function.\\
	This more general formulation of marginal models with stochastic regressors allows for a genuine extension beyond the case of independent individuals. Indeed, constant regressors imply that the conditional variances are also constant, which points to independence in the ``time index'' $i \ge 1.$\\
	Commonly used link functions include
	\begin{enumerate}
		\item $\mu(x)=x,$ for linear regression,
		\item $\mu(x)=\exp(x),$ for log regression for count data,
		\item $\mu(x) = \Phi(x),$ for probit regression for binary data.
	\end{enumerate} 
	We use the notation ${\mu}_i(\boldsymbol{\beta})=(\mu_{i1}(\boldsymbol{\beta}), \ldots, \mu_{im_i}(\boldsymbol{\beta}))^T$ and $\textbf{A}_i(\boldsymbol{\beta})$ for the diagonal (random) matrix with elements $\sigma_{i1}^2(\boldsymbol{\beta}), \ldots, \sigma_{im_i}^2(\boldsymbol{\beta})$, $i \ge 1$. Let $\Sigma_i^{(c)}(\boldsymbol{\beta})={\rm E}_{\boldsymbol{\beta}}[({\bf y}_i-\mu_i(\boldsymbol{\beta}))({\bf y}_i-\mu_i(\boldsymbol{\beta}))^T|\mathcal{F}_{i-1}]$, and $\overline{\textbf{R}}_i^{(c)}(\boldsymbol{\beta})=\textbf{A}_i(\boldsymbol{\beta})^{-1/2}\Sigma_i^{(c)}(\boldsymbol{\beta})\textbf{A}_i(\boldsymbol{\beta})^{-1/2}$. When $\boldsymbol{\beta}=\boldsymbol{\beta}_0$, $\Sigma_i^{(c)}=\Sigma_i^{(c)}(\boldsymbol{\beta}_0)$ is the true (unknown) covariance matrix, and $\overline{\textbf{R}}_i^{(c)}=\overline{\textbf{R}}_i^{(c)}(\boldsymbol{\beta}_0)$ is the true (unknown) correlation matrix, $i \ge 1$.
	
	We note that, in contrast to Balan et al. (2010), our correlation matrices may have random entries on the diagonal. This constitutes a unified and more general view of correlation matrices and is made possible by allowing for randomness in the regressors.   
	
	For any $\boldsymbol{\beta},$ we assume that the sequence of true correlation matrices is asymptotically a.s. positive definite on $(\Omega, \cF, P_{\boldsymbol{\beta}})$, i.e.\\
	
	$(H') \ \textrm{there exist constants} \ 0<C(\boldsymbol{\beta}) \le m, \ \textrm{such that} $
	\begin{eqnarray*}
		C(\boldsymbol{\beta}) \le \lambda_{\rm min}(\overline{{\bf R}}_i^{(c)}(\boldsymbol{\beta})) \le \lambda_{\rm max}(\overline{{\bf R}}_i^{(c)}(\boldsymbol{\beta})) \le m, \mbox{ a.s. } P_{\boldsymbol{\beta}},
	\end{eqnarray*}
	where the last inequality holds because $\overline{{\bf R}}_i^{(c)}(\boldsymbol{\beta})$ is a correlation matrix.
	
	As in Balan et al. (2010) we define the class of estimating functions 
	$$\mathcal{H}_n = \{ {\bf q}_n(\boldsymbol{\beta}) = \sum_{i=1}^n {\bf C}_i(\boldsymbol{\beta}) ({\bf y}_i - \mu_i(\boldsymbol{\beta})), \boldsymbol{\beta} \in \mathcal{T}\} $$
	where ${\bf C}_i(\boldsymbol{\beta})$ is a $p \times m_i$ random matrix with $\mathcal{F}_{i-1}$ measurable and continuously differentiable (in $\boldsymbol{\beta}$) entries. Furthermore, the following conditions need to be satisfied
	\begin{eqnarray}
	&& E_{\boldsymbol{\beta}} (|c_i^{jk}(\boldsymbol{\beta})|) < \infty, \ E_{\boldsymbol{\beta}} \left[\frac{\partial c_i^{jk}(\boldsymbol{\beta})}{\partial \boldsymbol{\beta}_l}(y_{ij} - \mu_{ij}(\boldsymbol{\beta}))\right] < \infty,  \nonumber   \\
	&& E_{\boldsymbol{\beta}}[c_i^{jk}(\boldsymbol{\beta})c_i^{rl}(\boldsymbol{\beta}) v^{kr}_i(\boldsymbol{\beta})]<\infty, \label{condition_new}
	\end{eqnarray}
	for any $\boldsymbol{\beta} \in \mathcal{T},$ $1 < k,r < m_i,$ $1 < j, l < p,$ where $c_i^{jk}(\boldsymbol{\beta})$ denotes the $(j,k)$ element of ${\bf C}_i(\boldsymbol{\beta})$ and $v^{kr}_i(\boldsymbol{\beta})$ is the $(k, r)$ component of $\Sigma_i^{(c)}(\boldsymbol{\beta}).$

	Within this class we are interested in estimating functions of the form
	\begin{equation}
	\textbf{g}_n^*(\boldsymbol{\beta})=\sum_{i=1}^n \textbf{X}_i^T \textbf{A}_i(\boldsymbol{\beta})^{1/2}\mathcal{R}_{i-1}^*(\boldsymbol{\beta})^{-1}\textbf{A}_i(\boldsymbol{\beta})^{-1/2} (\textbf{y}_i-\mu_i(\boldsymbol{\beta})), \label{GEE*}
	\end{equation}
	\noindent where $\{\mathcal{R}_{n}^*(\boldsymbol{\beta})\}_n$ is a sequence of $m_n \times m_n$ symmetric random matrices satisfying:\\
	\noindent (A) the matrix $\mathcal{R}_n^*(\boldsymbol{\beta})$ is positive definite for all $n \ge 1$ and $\boldsymbol{\beta} \in \mathcal{T}$
	
	\noindent (B) the entries of $\mathcal{R}_n^*(\boldsymbol{\beta})$ are $\mathcal{F}_n$ - measurable and continuously differentiable with respect to $\boldsymbol{\beta}$, for all $n \ge 1$ and $\boldsymbol{\beta} \in \mathcal{T}$.
	
	The estimating function 
	\begin{equation}
	\bar{g}_n(\boldsymbol{\beta})=\sum_{i=1}^n \textbf{X}_i^T \textbf{A}_i(\boldsymbol{\beta})^{1/2}\bar{R}_{i-1}^{(c)}(\boldsymbol{\beta})^{-1}\textbf{A}_i(\boldsymbol{\beta})^{-1/2} (\textbf{y}_i-\mu_i(\boldsymbol{\beta})), \label{quasiscore}
	\end{equation}
	plays an important role in the class $\mathcal{H}_n$; it is optimal (quasi-score), as shown in Proposition 3.6 of Balan et al. (2010). In Section 2 we describe sequences $\{{\bf g}_n^*(\boldsymbol{\beta})\}_{n\ge 1}$ which share, asymptotically the optimal properties of the score function \eqref{quasiscore}, i.e. are asymptotic score functions. This extends the results of Theorem 3.9 in Balan et al. (2010) to the case of stochastic regressors.
	
	A GEE estimator of the regression parameter $\boldsymbol{\beta}_0$ is an implicit solution of the equation ${\bf g}_n^*(\boldsymbol{\beta})=0$. 
	The matrices $\mathcal{R}_i^*(\boldsymbol{\beta})$ serve as proxies for the correlation matrices $\overline{{\bf R}}_i^{(c)}(\boldsymbol{\beta})$, $\boldsymbol{\beta} \in \mathcal{T}$, $i \ge 1$. They could be selected by the analyst, e.g. ``working independence'' estimating equations, where $\mathcal{R}_{i-1}^*(\boldsymbol{\beta}) = {\bf I}_{m_{i-1} \times m_{i-1}}$ is the identity matrix
	\begin{equation}
	{\bf g}_n^{\rm indep}(\boldsymbol{\beta}) = \sum_{i=1}^n {\bf X}_i^T (y_i - \mu_i(\boldsymbol{\beta})),  \nonumber
	\end{equation}
	or could be estimators of the true correlations, based on survey data, in which cases \eqref{GEE*} defines a pseudo-likelihood, 
	with $$\displaystyle{\mathcal{R}_n^*(\boldsymbol{\beta})=\frac{1}{n}\sum_{i=1}^n {\bf A}_i(\boldsymbol{ \boldsymbol{\beta}})^{-1/2}(y_i - \mu_i( \boldsymbol{\beta}))(y_i - \mu_i(\boldsymbol{\beta}))^T {\bf A}_i(\boldsymbol{\beta})^{-1/2}}.$$
	
	In all cases, and as stipulated in Liang and Zeger (1986), our theoretical results confirm that the accuracy to which the proxies represent the correlation matrices is irrelevant for the existence and strong consistency of the GEE estimators (see the hypotheses of our Theorem \ref{consistency}).  
	
	We assume throughout that the entries of $\mathcal{R}_i^*(\boldsymbol{\beta})^{-1}[y_{ij}- \mu_{ij}(\boldsymbol{\beta})]$ are $P_{\boldsymbol{\beta}}$-integrable, $\forall \boldsymbol{\beta} \in \mathcal{T}$, $j=1, \ldots m_i$, $i \ge 1$. Our model assumptions imply that the residuals form a martingale difference sequence. Combined with the measurability properties of the regressors, definition \eqref{GEE*} ensures that the sequence $({\bf g}^*_n)_{n \ge 1}$ is a transform martingale.
	
	For a $p \times p$ matrix, ${\bf A}$, we denote by $\displaystyle{\| {\bf A} \|=\sup_{\|{\bf x}\|=1}\|{\bf A}{\bf x}\|}$ its spectral norm and by $\displaystyle{\|| {\bf A} \||=\sup_{\|{\bf x}\|=1}|{\bf x}^T{\bf A}{\bf x}|}$ its numerical radius. We use the following inequality
	\begin{eqnarray*}
		&& \|| {\rm E} ({\bf A}) \|| = \sup_{\|{\bf x}\|=1} |{\bf x}^T {\rm E}({\bf A}) {\bf x}| \\
		&& \le \sup_{\|{\bf x}\|=1} {\rm E}|{\bf x}^T {\bf A} {\bf x}|  \\
		&& \le \sup_{\|{\bf x}\|=1} {\rm E} \sup_{\|{\bf x}\|=1}|{\bf x}^T {\bf A} {\bf x}| \\
		&& \le \sup_{\|{\bf x}\|=1} {\rm E} \|| {\bf A} \||.
	\end{eqnarray*}
	In addition, since $\||{\bf A} \|| \le \|{\bf A} \| \le  2 \|| {\bf A}\||$ we have $\| {\rm E} ({\bf A})\| \le 2 {\rm E}\|{\bf A} \|$.
	
	
	
	\section{Optimal Estimating Equations}
	\label{optim}
	
	We first introduce some notations. 
	
	We define, for $\boldsymbol{\beta} \in \mathcal{T}$ 
	\begin{eqnarray}
	\textbf{H}_n^{ind}(\boldsymbol{\beta})       &:=& \sum_{i=1}^n {\rm E}_{\boldsymbol{\beta}}\left[ \textbf{X}_i^T {\textbf{A}_i}(\boldsymbol{\beta}) \textbf{X}_i\right], \ n \ge 1  \nonumber \\
	\textbf{H}_n^*(\boldsymbol{\beta})           &:=& -\textrm{E}_{\boldsymbol{\beta}}\left[\frac{\partial\textbf{g}_n^*(\boldsymbol{\beta})}{\partial \boldsymbol{\beta}^T}\right]
	= \sum_{i=1}^n {\rm E}_{\boldsymbol{\beta}}\left[ \textbf{X}_i^T {\textbf{A}_i}(\boldsymbol{\beta})^{1/2} \mathcal{R}_{i-1}^*(\boldsymbol{\beta})^{-1} {\textbf{A}_i}(\boldsymbol{\beta})^{1/2} \textbf{X}_i\right] \nonumber \\
	&:=& \sum_{i=1}^n {\bf K}_i^*(\boldsymbol{\beta}), \ n \ge 1 \label{Hn^steluta} \\
	{\bf H}^*_{n_0, n}(\boldsymbol{\beta})       &:=& {\bf H}_n^*(\boldsymbol{\beta}) - {\bf H}^*_{n_0-1}(\boldsymbol{\beta}), \ 1 < n_0 < n \nonumber \\                          
	\overline{\textbf{M}}_n(\boldsymbol{\beta})  &:=& \textrm{Cov}_{\boldsymbol{\beta}}\left[\overline{\textbf{g}}_n(\boldsymbol{\beta})  \right] 
	=\sum_{i=1}^n {\rm E}_{\boldsymbol{\beta}}\left[\textbf{X}_i^T {\textbf{A}_i}(\boldsymbol{\beta})^{1/2} \overline{\textbf{R}}_{i}^{(c)}(\boldsymbol{\beta})^{-1}	  {\textbf{A}_i}(\boldsymbol{\beta})^{1/2} \textbf{X}_i \right] \nonumber \\
	&:=& \sum_{i=1}^n \overline{{\bf L}}_i(\boldsymbol{\beta}), \ n \ge 1 \label{Mn-bara} \\
	\overline{\bf M}_{n_0, n}(\boldsymbol{\beta})&:=& \overline{\bf M}_n(\boldsymbol{\beta}) - \overline{\bf M}_{n_0-1}(\boldsymbol{\beta}), \ 1 < n_0 < n \nonumber \\                          
	\textbf{M}_n^*(\boldsymbol{\beta})           &:=& \textrm{Cov}_{\boldsymbol{\beta}}\left[\textbf{g}_n^*(\boldsymbol{\beta})  \right] \nonumber \\
	&=& \sum_{i=1}^n {\rm E}_{\boldsymbol{\beta}}\left[ \textbf{X}_i^T {\textbf{A}_i}(\boldsymbol{\beta})^{1/2} \mathcal{R}_{i-1}^*(\boldsymbol{\beta})^{-1} \overline{\textbf{R}}_{i}^{(c)}(\boldsymbol{\beta})  \mathcal{R}_{i-1}^*(\boldsymbol{\beta})^{-1} {\textbf{A}_i}(\boldsymbol{\beta})^{1/2} \textbf{X}_i\right], \nonumber 
	\end{eqnarray}
	Let $(\boldsymbol{\delta}_i)_{i \ge 1}$ be a sequence of $p \times m_i$ random matrices, $\|{\boldsymbol{\delta}}_i\| \le d,$ for some $d>0.$ We define ${\bf Y}_i(\boldsymbol{\beta}, \boldsymbol{\delta}_i)^T: = ({\bf X}_i + \boldsymbol{\delta}_i)^T {\bf A}_i(\boldsymbol{\beta}, \boldsymbol{\delta}_i)^{1/2},$ where ${\bf A}_i(\boldsymbol{\beta}, \boldsymbol{\delta}_i)$ is obtained from ${\bf A}_i(\boldsymbol{\beta})$ by substituting ${\bf x}_{ij} + \boldsymbol{\delta}_{ij}$ for ${\bf x}_{ij}$ in $\mu'({\bf x}_{ij}^T\boldsymbol{\beta}),$ $j=1, \ldots, m_i,$ and let ${\bf Y}_i(\boldsymbol{\beta}):={\bf Y}_i(\boldsymbol{\beta}, {\bf 0}),$ $i\ge 1.$ Similarly, when $\mathcal{R}_i^*(\boldsymbol{\beta})$ depends on ${\bf X}_j,$ $j=1, \ldots, i$, we replace ${\bf X}_j$ by ${\bf X}_j + \boldsymbol{\delta}_j$ to obtain $\mathcal{R}_i^*(\boldsymbol{\beta}, \boldsymbol{\delta}_i),$ $i \ge 1.$ Note that ${\bf R}_i^{(c)}(\boldsymbol{\beta})$ does not depend on ${\bf X},$ $i \ge 1.$ We define
	\begin{equation}
	{\bf g}_n^*(\boldsymbol{\beta}, \boldsymbol{\delta}) := \sum_{i=1}^n {\bf Y}_i^T(\boldsymbol{\beta}, \boldsymbol{\delta}_i) \mathcal{R}_{i-1}^*(\boldsymbol{\beta}, \boldsymbol{\delta}_i)^{-1} {\bf A}_i(\boldsymbol{\beta}, \boldsymbol{\delta}_i)^{-1/2} ({\bf y}_i - \mu_i(\boldsymbol{\beta})) \label{gn^steluta}
	\end{equation}
	We can now state the main result of this section.
	\begin{theorem}
		\label{theo_opt}
		Let $\{\mathcal{R}_n^*(\boldsymbol{\beta})\}_{n\ge1}$ be a sequence of $m_n \times m_n$ random matrices satisfying $(A)$ and $(B)$. Assume that $(H')$ holds and for all $\boldsymbol{\beta} \in \mathcal{T}$
		
		$(D) \quad \lambda_{\min}[{\textbf{H}}_{n}^{ind}(\boldsymbol{\beta})] \to \infty$
		
		$(R)$  there exists $K(\boldsymbol{\beta})>0$ such that $\inf \lambda_{\min}[\mathcal{R}_n^*(\boldsymbol{\beta})] \ge K(\boldsymbol{\beta})$ $P_{\boldsymbol{\beta}} $ a.s. 
		
		$(A_1) \quad \mathcal{R}_{n-1}^*(\boldsymbol{\beta})- \overline{\textbf{R}}_n^{(c)}(\boldsymbol{\beta}) \stackrel{P_{\boldsymbol{\beta}} }{\longrightarrow} 0 \quad \textrm{element-wise},$ as $n \to \infty$
		
		
		\indent  In addition, assume that there exists a sequence of $p \times m_i$ random matrices $\boldsymbol{\delta}_i,$ 
		such that
		
		$(A_2) \ \max\{\left\|\mathcal{R}_{i-1}^*(\boldsymbol{\beta}, \boldsymbol{\delta}_i)^{-1}- \mathcal{R}_{i-1}^*(\boldsymbol{\beta})^{-1}\right\|, \|{\bf Y}_i(\boldsymbol{\beta}, \boldsymbol{\delta}_i) - {\bf Y}_i(\boldsymbol{\beta})\| \} \le \frac{1}{2^i},$ $i \ge 1$
		
		$(A_3) \ \{\|{\bf Y}_i(\boldsymbol{\beta}, \boldsymbol{\delta}_i)\|^2\}_{i \ge 1} \mbox { is } {\rm E}_{\boldsymbol{\beta}} \mbox{- uniformly integrable },$ 
		
		$(A_4) \ \overline{{\bf L}}_i(\boldsymbol{\beta}, \boldsymbol{\delta}_i)^{-1} \mbox{ and } {\bf K}_i^*(\boldsymbol{\beta}, \boldsymbol{\delta}_i)^{-1}, \ i \ge 1 \mbox{ exist and}$  
		$$ \inf_{i}{\rm E}_{\boldsymbol{\beta}} \{ \lambda_{\min}[ {\bf Y}_i(\boldsymbol{\beta}, \boldsymbol{\delta}_i)^T {\bf Y}_i(\boldsymbol{\beta}, \boldsymbol{\delta}_i)]\} > 0.$$
		
		
		
		\noindent Then the sequences $\{\textbf{g}_n^*(\boldsymbol{\beta})\}_{n \ge 1}$ and $\{\textbf{g}_n^*(\boldsymbol{\beta}, \boldsymbol{\delta})\}_{n \ge 1}$ are asymptotic quasi-score sequences in $\{\mathcal{H}_n\}_{n \ge 1}$.
	\end{theorem}
	
	\begin{remark}
		{\rm Condition $(A_1)$ is condition $(C)$ in Balan et al. (2010). While $(A_1)$ is essential in proving the optimality of $\{{\bf g}_n^*\}_{n \ge 1},$ it is not required in our Theorem \ref{consistency}, where we prove the existence and strong consistency of a sequence of estimators $\{\hat{\boldsymbol{\beta}}_n\}_{n \ge 1}.$} 
	\end{remark}
	
	\begin{corollary}
		\label{corl}
		Assume that $(H'),$ $(D),$ $(R)$ and $(A_1)$ hold. Assume further that ${\mathcal{R}}_i^*(\boldsymbol{\beta})$ is continuous in ${\bf X}_j,$ $1 \le j \le i,$ ${\bf A}_i(\boldsymbol{\beta})$ is continuous in ${\bf X}_i,$ $i \ge 1$ and that $(A_3')$ and $(A_4')$ hold for any $\boldsymbol{\beta} \in \mathcal{T},$
		
		$(A_3') \ \{\|{\bf Y}_i(\boldsymbol{\beta})\|^2\}_{i \ge 1} \mbox { is } {\rm E}_{\boldsymbol{\beta}} \mbox{- uniformly integrable },$ 
		
		$(A_4') \ \inf_{i}{\rm E}_{\boldsymbol{\beta}} \{ \lambda_{\min}[ {\bf Y}_i(\boldsymbol{\beta})^T {\bf Y}_i(\boldsymbol{\beta})]\} > 0.$
		
		Then there exists a sequence of matrices $\{\boldsymbol{\delta}_i\}_{i \ge 1}$ that satisfy the conditions of Theorem \ref{theo_opt} and therefore both $\{{\bf g}_n^*(\boldsymbol{\beta})\}_{n \ge 1}$ and $\{{\bf g}_n^*(\boldsymbol{\beta}, \boldsymbol{\delta})\}_{n \ge 1}$ are asymptotic quasi-score sequences in $\{\mathcal{H}_n\}_{n \ge 1}.$
	\end{corollary}
	
	{\bf Proof}. The continuity assumptions imply that there exists a sequence of $p \times m_i$ random matrices $\{\boldsymbol{\delta}_i\}_{i \ge 1}$ such that $(A_2)$ holds. These matrices can be chosen such that, using $(A_4'),$ $(A_4)$ holds. By $(A_2),$ the uniform integrability conditions $(A_3)$ and $(A_4')$ are equivalent. Thus, the hypothesis of Corollary \ref{corl} imply the conditions spelled out in Theorem \ref{theo_opt} and therefore its conclusion. $\Box$
	
	\vspace{3mm}
	
	\begin{remark}
		{\rm It is always possible to find matrices $\{\boldsymbol{\delta}_i\}_{i \ge 1}$ such that the inverse matrices in $(A_4)$ exist. However, for $\boldsymbol{\delta} \equiv {\bf 0},$ these conditions should be imposed in Theorem 3.9 of Balan et al. (2010). The continuity assumption in Corollary \ref{corl} not only ensures that the matrices $(A_4)$ exist, it also allows us to define a large class of sequences that are asymptotic quasi-scores in $\mathcal{H}_n,$ $n \ge 1,$ as long as the norm of $\boldsymbol{\delta}_i$ is ``small'', $i \ge 1.$ This means that even with slightly misspecified regressors ${\bf X}_i,$ $i \ge 1$ in the model, we can still obtain efficient estimators of the regression parameter.}
	\end{remark}
	
	{\bf Proof of Theorem \ref{theo_opt}}. As in \cite{balan-dumitrescu-schiopu10}, since $\overline {\textbf{g}}_n$ is an (asymptotic) quasi-score in $\mathcal{H}_n$, by Proposition 5.5 of \cite{heyde97} it is enough to show that
	\begin{eqnarray}
	&&\frac{\det \textbf{H}_n^*(\boldsymbol{\beta})} {\det \overline{\textbf{M}}_n(\boldsymbol{\beta})} \to 1, \quad \frac{\det \textbf{H}_n^*(\boldsymbol{\beta}, \boldsymbol{\delta}_i)} {\det \overline{\textbf{M}}_n(\boldsymbol{\beta}, \boldsymbol{\delta}_i)} \to 1 \quad \forall \boldsymbol{\beta} \in \mathcal{T}, \label{eq3} \\
	&&\frac{\det \textbf{M}_n^*(\boldsymbol{\beta})} {\det \overline{\textbf{M}}_n(\boldsymbol{\beta})} \to 1, \quad \frac{\det \textbf{M}_n^*(\boldsymbol{\beta}, \boldsymbol{\delta}_i)} {\det \overline{\textbf{M}}_n(\boldsymbol{\beta}, \boldsymbol{\delta}_i)} \to 1 \quad \forall \boldsymbol{\beta} \in \mathcal{T}. \label{eq4}
	\end{eqnarray}
	We proceed with the proof of \eqref{eq3}. The proof of \eqref{eq4} is similar and is omitted. 
	From $(A_4)$ we have,
	\begin{equation}
	\max \{\lambda_{\max}[\overline{{\bf L}}_i(\boldsymbol{\beta},\boldsymbol{\delta}_i)^{-1}], \lambda_{\max}[{\bf K}_i^*(\boldsymbol{\beta}, \boldsymbol{\delta}_i)^{-1}] \} < \infty \label{eq8}
	\end{equation} 
	First, we show that 
	\begin{equation}
	\sup_{i}\|\overline{{\bf L}}_i(\boldsymbol{\beta}, \boldsymbol{\delta}_i)^{-1}\| = \sup_{i}\lambda_{\max}[\overline{{\bf L}}_i(\boldsymbol{\beta}, \boldsymbol{\delta}_i)^{-1}]< \infty, \label{norm.inv.L.bounded}
	\end{equation}
	or $\inf_{i}\lambda_{\min}[\overline{{\bf L}}_i(\boldsymbol{\beta}, \boldsymbol{\delta}_i)] > 0.$ Indeed, let ${\bf x}$ be a vector of norm 1 such that ${\bf x}^T \overline{{\bf L}}_i(\boldsymbol{\beta}, \boldsymbol{\delta}_i){\bf x} =\lambda_{\min}[\overline{{\bf L}}_i(\boldsymbol{\beta}, \boldsymbol{\delta}_i)^{-1}].$ Since the absolute value of each entry of $\overline{{\bf R}}_i^{(c)}(\boldsymbol{\beta})$ is less than 1,
	\begin{eqnarray*}
		\lambda_{\min}[\overline{{\bf L}}_i(\boldsymbol{\beta},\boldsymbol{\delta}_i)] &=& {\rm E}_{\boldsymbol{\beta}}[{\bf x}^T {\bf Y}_i(\boldsymbol{\beta}, \boldsymbol{\delta}_i)^T \overline{\textbf{R}}_{i}^{(c)}(\boldsymbol{\beta})^{-1} {\bf Y}_i(\boldsymbol{\beta}, \boldsymbol{\delta}_i) {\bf x}] \\
		&\ge& {\rm E}_{\boldsymbol{\beta}}[\lambda_{\min}(\overline{\textbf{R}}_{i}^{(c)}(\boldsymbol{\beta})^{-1}) {\bf x}^T {\bf Y}_i(\boldsymbol{\beta}, \boldsymbol{\delta}_i)^T {\bf Y}_i(\boldsymbol{\beta}, \boldsymbol{\delta}_i){\bf x}] \\
		&=& {\rm E}_{\boldsymbol{\beta}}\left[\frac{1}{\lambda_{\max}(\overline{\textbf{R}}_{i}^{(c)}(\boldsymbol{\beta}))} {\bf x}^T {\bf Y}_i(\boldsymbol{\beta}, \boldsymbol{\delta}_i)^T {\bf Y}_i(\boldsymbol{\beta}, \boldsymbol{\delta}_i){\bf x} \right] \\
		&\ge& \frac{1}{m} {\rm E}_{\boldsymbol{\beta}} \{ \lambda_{\min}[ {\bf Y}_i(\boldsymbol{\beta}, \boldsymbol{\delta}_i)^T {\bf Y}_i(\boldsymbol{\beta}, \boldsymbol{\delta}_i)]\}.
	\end{eqnarray*}
	Thus, $\inf_{i} \lambda_{\min}[\overline{{\bf L}}_i(\boldsymbol{\beta},\boldsymbol{\delta}_i)]>0 $ by $(A_4).$
	
	
	For any $p \times 1$ vector ${\bf x}$, we have the following inequalities
	\begin{eqnarray}
	&&\min_{n_0 \le i \le n} \lambda_{\min}[\overline{{\bf L}}_i(\boldsymbol{\beta},\boldsymbol{\delta}_i)^{-1}{\bf K}_i^*(\boldsymbol{\beta},\boldsymbol{\delta}_i)] {\bf x}^T\overline{\textbf{M}}_{n_0,n}(\boldsymbol{\beta},\boldsymbol{\delta}){\bf x} \le {\bf x}^T\textbf{H}_{n_0,n}^*(\boldsymbol{\beta},\boldsymbol{\delta}){\bf x} \nonumber\\
	&&\le \max_{n_0 \le i \le n} \lambda_{\max}[\overline{{\bf L}}_i(\boldsymbol{\beta},\boldsymbol{\delta}_i)^{-1}{\bf K}_i^*(\boldsymbol{\beta},\boldsymbol{\delta}_i)] {\bf x}^T\overline{\textbf{M}}_{n_0,n}(\boldsymbol{\beta},\boldsymbol{\delta}){\bf x}, \label{eval_M_H}
	\end{eqnarray}
	where $\overline{\textbf{M}}_{n_0,n}(\boldsymbol{\beta},\boldsymbol{\delta}) := \sum_{i=n_0}^n \overline{{\bf L}}_i(\boldsymbol{\beta},\boldsymbol{\delta}_i)$ and $\textbf{H}_{n_0,n}^*(\boldsymbol{\beta},\boldsymbol{\delta}) := \sum_{i=n_0}^n {\bf K}_i^*(\boldsymbol{\beta},\boldsymbol{\delta}_i).$
	Next, we prove that as $ i \to \infty$,
	\begin{equation}
	\overline{{\bf L}}_i(\boldsymbol{\beta},\boldsymbol{\delta}_i)^{-1}{\bf K}_i^*(\boldsymbol{\beta},\boldsymbol{\delta}_i) \to {\bf I}. \label{convergenceLK}
	\end{equation}
	Combining \eqref{eval_M_H} and \eqref{convergenceLK}, we obtain that, for any $\varepsilon>0$, there exists $n_0= n_0(\varepsilon, \boldsymbol{\beta})$, such that 
	for any $n_m \ge n_0$ and $n > n_m,$
	\begin{equation}
	(1-\varepsilon) {\bf x}^T\overline{\textbf{M}}_{n_0,n}(\boldsymbol{\beta}, \boldsymbol{\delta}){\bf x} \le {\bf x}^T\textbf{H}_{n_0,n}^*(\boldsymbol{\beta}, \boldsymbol{\delta}){\bf x}  \le (1+\varepsilon) {\bf x}^T\overline{\textbf{M}}_{n_0,n}(\boldsymbol{\beta}, \boldsymbol{\delta}){\bf x}, \label{noua_eq}
	\end{equation}
	from which we derive
	\begin{equation}
	(1-\varepsilon)^p \le \frac{\det \textbf{H}_{n_0,n}^*(\boldsymbol{\beta}, \boldsymbol{\delta})} {\det \overline{\textbf{M}}_{n_0,n}(\boldsymbol{\beta}, \boldsymbol{\delta})} \le (1+\varepsilon)^p, \ \mbox{for any } n > n_m. \label{evaluation_partial_sum}
	\end{equation}
	We proceed with the proof of \eqref{convergenceLK}. We show that, as $i \to \infty$
	\begin{equation}
	\|\overline{{\bf L}}_i(\boldsymbol{\beta},\boldsymbol{\delta}_i)^{-1}{\bf K}_i^*(\boldsymbol{\beta},\boldsymbol{\delta}_i) - {\bf I}\| \le \|\overline{{\bf L}}_i(\boldsymbol{\beta},\boldsymbol{\delta}_i)^{-1}\| \cdot \|{\bf K}_i^*(\boldsymbol{\beta},\boldsymbol{\delta}_i)-\overline{{\bf L}}_i(\boldsymbol{\beta},\boldsymbol{\delta}_i)\| \to 0. \nonumber
	\end{equation}
	The first factor on the right hand side, is bounded by \eqref{norm.inv.L.bounded}. To complete the proof of \eqref{convergenceLK}, we now prove that $\|{\bf K}_i^*(\boldsymbol{\beta},\boldsymbol{\delta}_i)-\overline{{\bf L}}_i(\boldsymbol{\beta},\boldsymbol{\delta}_i)\| \to 0,$ $i \to \infty$.
	We have
	\begin{eqnarray*}
		&& \left\|{\rm E}_{\boldsymbol{\beta}}\left\{ {\bf Y}_i(\boldsymbol{\beta}, \boldsymbol{\delta}_i)^T \left[\mathcal{R}_{i-1}^*(\boldsymbol{\beta}, \boldsymbol{\delta}_i)^{-1}- \mathcal{R}_{i-1}^*(\boldsymbol{\beta})^{-1} + \mathcal{R}_{i-1}^*(\boldsymbol{\beta})^{-1} - \overline{\textbf{R}}_{i}^{(c)}(\boldsymbol{\beta})^{-1}\right]{\bf Y}_i(\boldsymbol{\beta}, \boldsymbol{\delta}_i)\right\} \right\|  \\
		&& \le \left\|{\rm E}_{\boldsymbol{\beta}}\left\{ {\bf Y}_i(\boldsymbol{\beta}, \boldsymbol{\delta}_i)^T \left[\mathcal{R}_{i-1}^*(\boldsymbol{\beta}, \boldsymbol{\delta}_i)^{-1}- \mathcal{R}_{i-1}^*(\boldsymbol{\beta})^{-1}\right]{\bf Y}_i(\boldsymbol{\beta}, \boldsymbol{\delta}_i)\right\} \right\|  \\
		&& + \left\|{\rm E}_{\boldsymbol{\beta}}\left\{ {\bf Y}_i(\boldsymbol{\beta}, \boldsymbol{\delta}_i)^T \left[\mathcal{R}_{i-1}^*(\boldsymbol{\beta})^{-1}- \overline{\textbf{R}}_{i}^{(c)}(\boldsymbol{\beta})^{-1} \right]{\bf Y}_i(\boldsymbol{\beta}, \boldsymbol{\delta}_i)\right\} \right\|  \\
		&& := T_1(i) + T_2(i). 
	\end{eqnarray*}
	Using $(A_2)$, we obtain the inequalities 
	\begin{eqnarray*}
		&& T_1(i) \le 2{\rm E}_{\boldsymbol{\beta}} \left\{\left\| {\bf Y}_i(\boldsymbol{\beta}, \boldsymbol{\delta}_i)^T \left[\mathcal{R}_{i-1}^*(\boldsymbol{\beta}, \boldsymbol{\delta}_i)^{-1}- \mathcal{R}_{i-1}^*(\boldsymbol{\beta})^{-1}\right]{\bf Y}_i(\boldsymbol{\beta}, \boldsymbol{\delta}_i) \right\| \right\} \\
		&& \le 2{\rm E}_{\boldsymbol{\beta}} \left\{\| {\bf Y}_i(\boldsymbol{\beta}, \boldsymbol{\delta}_i) \|^2 \left\|\mathcal{R}_{i-1}^*(\boldsymbol{\beta}, \boldsymbol{\delta}_i)^{-1}- \mathcal{R}_{i-1}^*(\boldsymbol{\beta})^{-1}\right\| \right\} \le 2\frac{1}{2^i} \sup_{i} {\rm E}_{\boldsymbol{\beta}} [\|{\bf Y}_i(\boldsymbol{\beta}, \boldsymbol{\delta}_i) \|^2].
	\end{eqnarray*} 
	
	By $(A_3),$ we obtain $T_1(i) \to 0,$ as $i \to \infty.$ On the other hand,  
	\begin{eqnarray*}
		&& T_2(i) \le  2{\rm E}_{\boldsymbol{\beta}}\{\| {\bf Y}_i(\boldsymbol{\beta}, \boldsymbol{\delta}_i)^T [\mathcal{R}_{i-1}^*(\boldsymbol{\beta})^{-1}- \overline{\textbf{R}}_{i}^{(c)}(\boldsymbol{\beta})^{-1} ]{\bf Y}_i(\boldsymbol{\beta}, \boldsymbol{\delta}_i) \|\} \\
		&& \le 2{\rm E}_{\boldsymbol{\beta}}\{\| {\bf Y}_i(\boldsymbol{\beta}, \boldsymbol{\delta}_i)\|^2 \|\mathcal{R}_{i-1}^*(\boldsymbol{\beta})^{-1}- \overline{\textbf{R}}_{i}^{(c)}(\boldsymbol{\beta})^{-1} \|\} \\
		&& \le 2{\rm E}_{\boldsymbol{\beta}}\{\| {\bf Y}_i(\boldsymbol{\beta}, \boldsymbol{\delta}_i)\|^2 \|\mathcal{R}_{i-1}^*(\boldsymbol{\beta})^{-1}- \overline{\textbf{R}}_{i}^{(c)}(\boldsymbol{\beta})^{-1} \| {\bf 1}_{\{\|\mathcal{R}_{i-1}^*(\boldsymbol{\beta})^{-1}- \overline{\textbf{R}}_{i}^{(c)}(\boldsymbol{\beta})^{-1}\| < \varepsilon\} }\} \\
		&& +  2 {\rm E}_{\boldsymbol{\beta}}\{\| {\bf Y}_i(\boldsymbol{\beta}, \boldsymbol{\delta}_i)\|^2 \|\mathcal{R}_{i-1}^*(\boldsymbol{\beta})^{-1}- \overline{\textbf{R}}_{i}^{(c)}(\boldsymbol{\beta})^{-1} \| {\bf 1}_{\{\|\mathcal{R}_{i-1}^*(\boldsymbol{\beta})^{-1}- \overline{\textbf{R}}_{i}^{(c)}(\boldsymbol{\beta})^{-1}\| \ge \varepsilon\} }\}
	\end{eqnarray*}
	which gives
	\begin{eqnarray*}
		T_2(i) \le 2\varepsilon \sup_{i}{\rm E}_{\boldsymbol{\beta}}\{\| {\bf Y}_i(\boldsymbol{\beta}, \boldsymbol{\delta}_i)\|^2 \}+ 2M_1 \eta, \ i \ge i_0(\varepsilon).
	\end{eqnarray*}
	We obtained the bound for the second term on the right hand side by using conditions $(H'),$ $(R),$ $(A_1)$ and $(A_3).$ 
	This proves that $T_2(i) \to 0$ as $i \to \infty$ and concludes the proof of \eqref{convergenceLK}. Since all the eigenvalues of the matrix $\overline{{\bf L}}_i(\boldsymbol{\beta},\boldsymbol{\delta}_i)^{-1}{\bf K}_i^*(\boldsymbol{\beta},\boldsymbol{\delta}_i)$ converge to 1, for any $\varepsilon > 0,$ there exists $n_0$ such that 
	\begin{equation}
	1- \varepsilon \le \min_{n_0 \le i \le n} \lambda_{\min}[\overline{{\bf L}}_i(\boldsymbol{\beta},\boldsymbol{\delta}_i)^{-1}{\bf K}_i^*(\boldsymbol{\beta},\boldsymbol{\delta}_i)] \le \max_{n_0 \le i \le n} \lambda_{\max}[\overline{{\bf L}}_i(\boldsymbol{\beta},\boldsymbol{\delta}_i)^{-1}{\bf K}_i^*(\boldsymbol{\beta},\boldsymbol{\delta}_i)] \le 1+\varepsilon, \nonumber
	\end{equation}
	Combining these inequalities with \eqref{eval_M_H} we obtain \eqref{noua_eq}, and thus \eqref{evaluation_partial_sum}. 
	Actually, our first goal is to obtain inequalities similar to those in \eqref{evaluation_partial_sum}, when all $\boldsymbol{\delta}_i=0,$ $i \ge 1.$ We proceed as follows.
	
	For any $p \times 1$ vector ${\bf x}$, we have
	\begin{eqnarray}
	&&\min_{n_0 \le i \le n} \lambda_{\min}[{\bf K}_i^*(\boldsymbol{\beta}, \boldsymbol{\delta}_i)^{-1}{\bf K}_i^*(\boldsymbol{\beta})] {\bf x}^T\textbf{H}_{n_0,n}^*(\boldsymbol{\beta}, \boldsymbol{\delta}){\bf x} \le {\bf x}^T\textbf{H}_{n_0,n}^*(\boldsymbol{\beta}){\bf x} \nonumber\\
	&&\le \max_{n_0 \le i \le n} \lambda_{\max}[{\bf K}_i^*(\boldsymbol{\beta}, \boldsymbol{\delta}_i)^{-1}{\bf K}_i^*(\boldsymbol{\beta})] {\bf x}^T\textbf{H}_{n_0,n}^*(\boldsymbol{\beta}, \boldsymbol{\delta}){\bf x}. \label{eval_H_Hd}
	\end{eqnarray}
	and
	\begin{eqnarray}
	&&\min_{n_0 \le i \le n} \lambda_{\min}[\overline{{\bf L}}_i(\boldsymbol{\beta}, \boldsymbol{\delta}_i)^{-1}\overline{{\bf L}}_i(\boldsymbol{\beta})] {\bf x}^T\overline{\textbf{M}}_{n_0,n}(\boldsymbol{\beta}, \boldsymbol{\delta}){\bf x} \le {\bf x}^T\overline{\textbf{M}}_{n_0,n}(\boldsymbol{\beta}){\bf x} \nonumber\\
	&&\le \max_{n_0 \le i \le n} \lambda_{\max}[\overline{{\bf L}}_i(\boldsymbol{\beta},\boldsymbol{\delta}_i)^{-1}\overline{{\bf L}}_i(\boldsymbol{\beta})] {\bf x}^T\overline{\textbf{M}}_{n_0,n}(\boldsymbol{\beta}, \boldsymbol{\delta}){\bf x}. \label{eval_M_Md}
	\end{eqnarray}
	We will prove that, as $ i \to \infty$,
	\begin{equation}
	{\bf K}_i^*(\boldsymbol{\beta}, \boldsymbol{\delta}_i)^{-1}{\bf K}_i^*(\boldsymbol{\beta}) \longrightarrow {\bf I}, \label{convergenceKKd}
	\end{equation}
	and hence all the eigenvalues of the matrix ${\bf K}_i^*(\boldsymbol{\beta}, \boldsymbol{\delta}_i)^{-1}{\bf K}_i^*(\boldsymbol{\beta})$ converge to 1. If \eqref{convergenceKKd} holds then, for $\varepsilon>0$, there exists $n_1$ such that 
	\begin{equation}
	1- \varepsilon \le \min_{n_1 \le i \le n} \lambda_{\min}[{\bf K}_i^*(\boldsymbol{\beta}, \boldsymbol{\delta}_i)^{-1}{\bf K}_i^*(\boldsymbol{\beta})] \le \max_{n_1 \le i \le n} \lambda_{\max}[{\bf K}_i^*(\boldsymbol{\beta}, \boldsymbol{\delta}_i)^{-1}{\bf K}_i^*(\boldsymbol{\beta})] \le 1+\varepsilon. \nonumber
	\end{equation}
	These inequalities, combined with \eqref{eval_H_Hd} imply first that for all $n_m,$ $n,$ $n_1 \le n_m< n$
	\begin{equation}
	\label{eq222}
	(1-\varepsilon) {\bf x}^T\textbf{H}_{n_m,n}^*(\boldsymbol{\beta}, \boldsymbol{\delta}){\bf x} \le {\bf x}^T\textbf{H}_{n_m,n}^*(\boldsymbol{\beta}){\bf x} \le (1+\varepsilon) {\bf x}^T\textbf{H}_{n_m,n}^*(\boldsymbol{\beta}, \boldsymbol{\delta}){\bf x}.
	\end{equation}
	and then 
	\begin{equation}
	\label{evaluation_partial_sum_H_Hd}
	(1-\varepsilon)^p  \le \frac{\det \textbf{H}_{n_m,n}^*(\boldsymbol{\beta})}{\det \textbf{H}_{n_m,n}^*(\boldsymbol{\beta}, \boldsymbol{\delta})} \le (1+\varepsilon)^p, \ \mbox{for any } n > n_m. 
	\end{equation}
	Similarly, if 
	\begin{equation}
	\overline{{\bf L}}_i(\boldsymbol{\beta},\boldsymbol{\delta}_i)^{-1}\overline{{\bf L}}_i(\boldsymbol{\beta}) \longrightarrow {\bf I},  \label{new_label2}
	\end{equation} 
	then from \eqref{eval_M_Md} there exists $n_2$ such that, for any $n_m \ge n_2$ and $n \ge n_m,$ 
	\begin{equation}
	(1-\varepsilon){\bf x}^T\overline{\textbf{M}}_{n_m,n}(\boldsymbol{\beta}, \boldsymbol{\delta}){\bf x} \le {\bf x}^T\overline{\textbf{M}}_{n_m,n}(\boldsymbol{\beta}){\bf x} \le (1+\varepsilon){\bf x}^T\overline{\textbf{M}}_{n_m,n}(\boldsymbol{\beta}, \boldsymbol{\delta}){\bf x}, \label{eq111}
	\end{equation}
	and 
	\begin{equation}
	\label{evaluation_partial_sum_M_Md}
	(1-\varepsilon)^p \le \frac{ \det \overline{\textbf{M}}_{n_m,n}(\boldsymbol{\beta})}{\det \overline{\textbf{M}}_{n_m,n}(\boldsymbol{\beta}, \boldsymbol{\delta})} \le (1+\varepsilon)^p , \ \mbox{for any } n > n_m. 
	\end{equation}
	Let
	\begin{equation}
	\frac{\det \textbf{H}_{n_m,n}^*(\boldsymbol{\beta})} {\det \overline{\textbf{M}}_{n_m,n}(\boldsymbol{\beta})} = \frac{\det \textbf{H}_{n_m,n}^*(\boldsymbol{\beta})} {\det \textbf{H}_{n_m,n}^*(\boldsymbol{\beta}, \boldsymbol{\delta})} \frac{\det \textbf{H}_{n_m,n}^*(\boldsymbol{\beta}, \boldsymbol{\delta})}{\det \overline{\textbf{M}}_{n_m,n}(\boldsymbol{\beta}, \boldsymbol{\delta})} \frac{\det \overline{\textbf{M}}_{n_m,n}(\boldsymbol{\beta}, \boldsymbol{\delta})}{\det \overline{\textbf{M}}_{n_m,n}(\boldsymbol{\beta})}. \nonumber
	\end{equation}
	Combining \eqref{evaluation_partial_sum}, \eqref{evaluation_partial_sum_H_Hd} and \eqref{evaluation_partial_sum_M_Md} we obtain 
	\begin{equation}
	(1-\varepsilon)^p(1-\varepsilon)^{p}(1+\varepsilon)^{-p} \le \frac{\det \textbf{H}_{n_m,n}^*(\boldsymbol{\beta})} {\det \overline{\textbf{M}}_{n_m,n}(\boldsymbol{\beta})} \le (1+\varepsilon)^p(1+\varepsilon)^{p}(1-\varepsilon)^{-p}, \ \mbox{for any } n > n_m. \nonumber
	\end{equation}
	where $n_3 \le n_m \le n,$ where $n_3 = \max\{n_0, n_1, n_2\}.$
	We now turn to the proof of \eqref{convergenceKKd} which closely follows the proof of \eqref{convergenceLK}. 
	The proof of \eqref{new_label2} is similar. We show that, as $i \to \infty$
	\begin{equation}
	\|{\bf K}_i^*(\boldsymbol{\beta},\boldsymbol{\delta}_i)^{-1}{\bf K}_i^*(\boldsymbol{\beta}) - {\bf I}\| \le \|{\bf K}_i^*(\boldsymbol{\beta},\boldsymbol{\delta}_i)^{-1}\| \cdot \|{\bf K}_i^*(\boldsymbol{\beta})-{\bf K}_i^*(\boldsymbol{\beta},\boldsymbol{\delta}_i)\| \to 0. \nonumber
	\end{equation}
	By $(A_2)$ and $(A_3)$ we have $ \sup_{i}{\rm E}_{\boldsymbol{\beta}}[ \|{\bf Y}_i(\boldsymbol{\beta})\|] < \infty.$\\
	
	A further use of $(A_2)$ and $(A_3)$ gives
	\begin{eqnarray*}
		&& \| {\bf K}_i^*(\boldsymbol{\beta})-{\bf K}_i^*(\boldsymbol{\beta},\boldsymbol{\delta}_i) \| \\
		&& = \| {\rm E}_{\boldsymbol{\beta}}[ {\bf Y}_i(\boldsymbol{\beta})^T \mathcal{R}_{i-1}^*(\boldsymbol{\beta})^{-1} {\bf Y}_i(\boldsymbol{\beta})] - {\rm E}_{\boldsymbol{\beta}}\left[ {\bf Y}_i(\boldsymbol{\beta}, \boldsymbol{\delta}_i)^T \mathcal{R}_{i-1}^*(\boldsymbol{\beta}, \boldsymbol{\delta}_i)^{-1} {\bf Y}_i(\boldsymbol{\beta}, \boldsymbol{\delta}_i)\right] \| \\
		&& \le 2{\rm E}_{\boldsymbol{\beta}} \|[{\bf Y}_i(\boldsymbol{\beta})^T - {\bf Y}_i(\boldsymbol{\beta}, \boldsymbol{\delta}_i)]^T\mathcal{R}_{i-1}^*(\boldsymbol{\beta})^{-1} {\bf Y}_i(\boldsymbol{\beta})\| \\
		&& +  2{\rm E}_{\boldsymbol{\beta}} \|{\bf Y}_i(\boldsymbol{\beta}, \boldsymbol{\delta}_i)^T [\mathcal{R}_{i-1}^*(\boldsymbol{\beta})^{-1}- \mathcal{R}_{i-1}^*(\boldsymbol{\beta}, \boldsymbol{\delta}_i)^{-1} ] {\bf Y}_i(\boldsymbol{\beta}, \boldsymbol{\delta}_i)\| \\
		&& + 2{\rm E}_{\boldsymbol{\beta}} \|{\bf Y}_i(\boldsymbol{\beta}, \boldsymbol{\delta}_i)^T \mathcal{R}_{i-1}^*(\boldsymbol{\beta})^{-1}[{\textbf{A}_i}(\boldsymbol{\beta})^{1/2} \textbf{X}_i -{\bf Y}_i(\boldsymbol{\beta}, \boldsymbol{\delta}_i) ] \| \\
		&\le& 2\frac{1}{2^i} K(\boldsymbol{\beta})^{-1} \sup_{i}{\rm E}_{\boldsymbol{\beta}}[ \|{\bf Y}_i(\boldsymbol{\beta})\|] + 2\frac{1}{2^i} \sup_i {\rm E}_{\boldsymbol{\beta}}[\|{\bf Y}_i(\boldsymbol{\beta}, \boldsymbol{\delta}_i)^2\|] + 4\frac{1}{2^i} K(\boldsymbol{\beta})^{-1} \sup_{i}{\rm E}_{\boldsymbol{\beta}}[ \|{\bf Y}_i(\boldsymbol{\beta}, \boldsymbol{\delta}_i)\|]\\
		&& \to 0, \ \mbox{as }i \to \infty,
	\end{eqnarray*}
	by (R) and because since all expectations are equibounded. 
	Thus \eqref{convergenceKKd} holds and the proof of \eqref{evaluation_partial_sum_H_Hd} is complete.\\ To complete the proof of the first part of \eqref{eq3}, we have to deal with the terms that are missing in \eqref{evaluation_partial_sum_H_Hd}.
	
	Since $\overline{{\bf M}}_{n}(\boldsymbol{\beta}) \ge m {\bf H}_{n}^{ind}(\boldsymbol{\beta}),$ $\lambda_{\min}[\overline{{\bf M}}_{n}(\boldsymbol{\beta})] \to \infty,$ as $n \to \infty$ by condition $(D),$ we obtain $\lambda_{\min}[\overline{{\bf M}}_{n_m, n}(\boldsymbol{\beta})] \to \infty,$ for any fixed $n_m$ as $n \to \infty.$ For a given $\varepsilon >0$ we can find $n_4 =n_4(\varepsilon,\boldsymbol{\beta}) > n_m,$ such that $\lambda_{\min}[\overline{{\bf M}}_{n_m, n}(\boldsymbol{\beta})] \ge \varepsilon^{-1} \lambda_{\max}[\overline{{\bf M}}_{n_m-1}(\boldsymbol{\beta})],$ for all $n \ge n_4.$ Thus, for all $n \ge n_4 > n_m.$
	
	\begin{equation}
	\nonumber
	\overline{{\bf M}}_{n_m, n}(\boldsymbol{\beta}) \le \overline{{\bf M}}_{n}(\boldsymbol{\beta}) \le (1+\varepsilon)\overline{{\bf M}}_{n_m, n}(\boldsymbol{\beta}) 
	\end{equation} 
	and
	\begin{equation}
	\label{eq_new}
	\det \overline{{\bf M}}_{n_m, n}(\boldsymbol{\beta})  \le  \det \overline{{\bf M}}_{n}(\boldsymbol{\beta})  \le  (1+\varepsilon)^p \det \overline{{\bf M}}_{n_m, n}(\boldsymbol{\beta}). 
	\end{equation}
	Combining \eqref{noua_eq}, \eqref{eq222} and \eqref{eq111} we obtain for all $n > n_m$
	\begin{eqnarray*}
		(1+\varepsilon)^{-1} (1-\varepsilon)^2 {\bf x}^T\overline{\textbf{M}}_{n_m, n}(\boldsymbol{\beta}){\bf x} \le {\bf x}^T\textbf{H}_{n_m, n}^*(\boldsymbol{\beta}){\bf x} \le (1+\varepsilon)^2 (1-\varepsilon)^{-1} {\bf x}^T\overline{\textbf{M}}_{n_m, n}(\boldsymbol{\beta}){\bf x}.
	\end{eqnarray*}
	These inequalities imply that $\lambda_{\min}[\textbf{H}_{n_m,n}^*(\boldsymbol{\beta})] \to \infty.$ Reasoning as above, there exists an integer $n_5 \ge  n_4,$ such that, for all $n \ge n_5 > n_m,$
	\begin{equation}
	\det {\bf H}_{n_m, n}^*(\boldsymbol{\beta}) \le \det {\bf H}_{n}^*(\boldsymbol{\beta}) \le (1 + \varepsilon)^p \det {\bf H}_{n_m, n}^*(\boldsymbol{\beta}), \label{eq25}
	\end{equation}
	Combining \eqref{eq_new} and \eqref{eq25} gives, for $n \ge n_5 > n_m$
	\begin{equation}
	\frac{1} {(1+\varepsilon)^p} \frac{\det {\bf H}_{n_m, n}^*(\boldsymbol{\beta})}{\det \overline{\textbf{M}}_{n_m,n}(\boldsymbol{\beta})}  \le \frac{\det {\bf H}_{n}^*(\boldsymbol{\beta})}{\det \overline{\textbf{M}}_{n_m,n}(\boldsymbol{\beta})} \le (1+ \varepsilon)^p \frac{\det {\bf H}_{n_m, n}^*(\boldsymbol{\beta})}{\det \overline{\textbf{M}}_{n_m,n}(\boldsymbol{\beta})}. \label{eq_eq}
	\end{equation}
	Finally, we obtain the first part of \eqref{eq3} from \eqref{evaluation_partial_sum_H_Hd} and \eqref{eq_eq}. To prove $\displaystyle{\frac{\det {\bf H}_{n}^*(\boldsymbol{\beta}, \boldsymbol{\delta})}{\det \overline{\bf M}_{n}(\boldsymbol{\beta}, \boldsymbol{\delta})}} \to 1,$ as $n \to \infty$ we proceed in a similar way. We note first that, by \eqref{eq222}, $\lambda_{\min}[{\bf H}_{m_n, n}^*(\boldsymbol{\beta})] \to \infty$ is equivalent to $\lambda_{\min}[{\bf H}_{m_n, n}^*(\boldsymbol{\beta}, \boldsymbol{\delta})] \to \infty,$ as $n \to \infty.$ Then $\overline{\bf M}_{m_n, n}(\boldsymbol{\beta})$ and ${\bf H}_{m_n, n}^*(\boldsymbol{\beta})$ can be replaced with $\overline{\bf M}_{m_n, n}(\boldsymbol{\beta}, \boldsymbol{\delta})$ and ${\bf H}_{m_n, n}^*(\boldsymbol{\beta}, \boldsymbol{\delta}),$ respectively in \eqref{eq_new} - \eqref{eq_eq}. Now the second part of \eqref{eq3} holds by \eqref{evaluation_partial_sum}.     $\Box$

	\section{Strong Consistency}
	\label{consist}
	
	From here on we work on a single probability space $(\Omega, \mathcal{F}, P_{\boldsymbol{\beta}_0})$ and all sigma-fields are assumed to be $P_{\boldsymbol{\beta}_0}$-complete. We consider $p$-dimensional random functions ${\bf q}_n(\boldsymbol{\beta})=\sum_{i=1}^n{\bf u}_i(\boldsymbol{\beta}),$ with $\boldsymbol{\beta} \in \cT$, which are continuously differentiable with respect to $\boldsymbol{\beta}$, ${n \geq 1}$. In keeping with an established convention, we omit writing $\boldsymbol{\beta}_0$ when it appears in the argument of a function, e.g., ${\bf q}_n = {\bf q}_n(\boldsymbol{\beta}_0)$, $n \ge 1$. We assume that each ${\bf q}_n$ is a square integrable martingale with $0$ mean and {\em a.s.} positive definite conditional variance, ${\bf V}_n=\sum_{i=1}^n {\rm E}({\bf u}_i{\bf u}_i^T |\cF_{i-1})$, $n \ge 1$. Let
	\begin{equation}
	\cD_n(\boldsymbol{\beta})=-\frac{\partial {\bf q}_n(\boldsymbol{\beta})}{\partial \boldsymbol{\beta}^T}, \quad  \boldsymbol{\beta} \in \cT. \nonumber
	\end{equation}
	First, we consider the simplest case of a longitudinal linear model with stochastic regressors as in Example 5.1 of Xie and Yang (2003).
	
	\begin{example}
		{\em Consider the estimating function ${\bf q}_n(\boldsymbol{\beta}) \in \mathcal{H}_n,$ $\beta \in \cT$
			$${\bf q}_n(\boldsymbol{\beta}) = \sum_{i=1}^n {\bf X}_i^T {\bf R}_{i}(\alpha)^{-1} ({\bf y}_i - {\bf X}_i \boldsymbol{\beta}),$$
			where ${\bf R}_i(\alpha)$ (which represents $\mathcal{R}_{i-1}^*(\boldsymbol{\beta})$) satisfies conditions $(A)$ and $(B),$ but does not depend on $\beta$ (see Xie and Yang, 2003). In our case, ${\bf x}_{ij}$ are random.
			The corresponding estimating equation can be solved explicitly and the solution is $$\widehat{\boldsymbol{\beta}}_n = (\sum_{i=1}^n {\bf X}_i^T {\bf R}_{i}(\alpha)^{-1} {\bf X}_i)^{-1} (\sum_{i=1}^n {\bf X}_i^T {\bf R}_{i}(\alpha)^{-1}{\bf y}_i).$$ 
			
			Then, $\widehat{\boldsymbol{\beta}}_n - \boldsymbol{\beta}_0 \to 0,$ a.s. is equivalent to $$(\sum_{i=1}^n {\bf X}_i^T {\bf R}_{i}(\alpha)^{-1}{\bf X}_i)^{-1} (\sum_{i=1}^n {\bf X}_i^T {\bf R}_{i}(\alpha)^{-1}\boldsymbol{\varepsilon}_i) \to 0 \mbox{ a.s.,}$$
			where $\boldsymbol{\varepsilon}_i= {\bf y}_i - {\bf X}_i \boldsymbol{\beta}_0$ and so the strong consistency of $\widehat{\boldsymbol{\beta}}_n $ is equivalent to 
			$${\bf H}_n^{-1} {\bf s}_n \to 0 \mbox{ a.s.,}$$
			where ${\bf H}_n = \sum_{i=1}^n {\bf X}_i^T {\bf R}_{i}(\alpha)^{-1}{\bf X}_i,$ and ${\bf s}_n= \sum_{i=1}^n {\bf X}_i^T {\bf R}_{i}(\alpha)^{-1}\boldsymbol{\varepsilon}_i \mbox{ is a martingale}.$ 
			Applying Theorem 4 of Lin (1994), with ${\bf F}_n = \sum_{i=1}^n {\bf X}_i^T {\bf R}_{i}(\alpha)^{-1}\overline{{\bf R}}_i^{(c)}(\boldsymbol{\beta}_0) {\bf R}_{i}(\alpha)^{-1}{\bf X}_i$ we obtain
			$${\bf F}_n^{-1} {\bf s}_n \to 0 \mbox{ a.s. and in } L_2,$$
			under the hypothesis $[\log \lambda_{\max}({\bf F}_n)]^{\nu} = o(\lambda_{\min}({\bf F}_n))$ for some $\nu>1,$ together with the assumption $E[{\bf X}_i^T {\bf R}_i(\alpha)^{-1}\boldsymbol{\Sigma}_i^{(c)} {\bf R}_i(\alpha)^{-1}{\bf X}_i] < \infty,$ for any $i.$ This last assumption guarantees square integrability. 
			
			If there exists a constant $c>0$ such that $ c \le \min_{i \le n} {\lambda}_{\min} [{\bf R}_{i}(\alpha)],$ then the almost sure convergence holds with the normalizer ${\bf H}_n$ in lieu of ${\bf F}_n.$
			
			We conclude that in this case strong consistency follows if $[\log \lambda_{\max}({\bf F}_n)]^{\nu} = o(\lambda_{\min}({\bf F}_n))$ holds for some $\nu >1.$ This result extends Corollary 3 of Lai and Wei (1982) to the longitudinal case.
		}
	\end{example}

	The following Lemma is an essential tool which allows us to obtain an a.s. convergence result for a general estimating functions, with random norming. 
	\begin{lemma}
		\label{random-norm-martingale-convergence}
		If $\lambda_{\min}({\bf V}_n) \to \infty$, then, for any $\displaystyle{\boldsymbol{\delta}>0}$
		\begin{equation}
		\frac{{\bf q}_n}{[\lambda_{\max}({\bf V}_n)]^{1/2+\boldsymbol{\delta}}} \longrightarrow 0, \mbox{a.s. and in }L_2. \nonumber
		\end{equation} 
	\end{lemma}
	{\bf Proof.} The proof is based on a result in Theorem 4 of Lin (1994). 
	
	For any $k \in \{1, \ldots p\}$ we show that  
	\begin{equation}
	\frac{{\bf q}_n^k}{[\lambda_{\max}({\bf V}_n)]^{1/2+\boldsymbol{\delta}}} \longrightarrow 0, \mbox{a.s. and in }L_2, \nonumber
	\end{equation} 
	where ${\bf q}_n^k$ is the $k$-th component of the vector ${\bf q}_n$. If ${\bf e}_k^T=(0, \ldots, 1, \ldots, 0)$ is the $p$-dimensional vector with the $k$-th component equal to 1, then ${\bf q}_n^k={\bf e}_k^T{\bf q}_n=\sum_{i=1}^n{\bf e}_k^T {\bf u}_i$ and ${\bf v}_n^k=\sum_{i=1}^n {\rm E}({\bf e}_k^T{\bf u}_i{\bf u}_i^T{\bf e}_k |\cF_{i-1})={\bf e}_k^T{\bf V}_n{\bf e}_k$. Since $\lambda_{\min}({\bf V}_n) \le {\bf v}_n^k$ 
	and $\lambda_{\min}({\bf V}_n) \to \infty$, we have ${\bf v}_n^k \to \infty$. This together with Theorem 4, (ii) in Lin (1994) (with $\alpha=1/2+\boldsymbol{\delta}$ and any $\nu>1/(2\boldsymbol{\delta})$) implies
	\begin{equation}
	\frac{{\bf q}_n^k}{({\bf v}_n^k)^{1/2+\boldsymbol{\delta}}} \longrightarrow 0, \mbox{a.s. and in }L_2. \nonumber
	\end{equation}
	The conclusion now follows, since ${\bf v}_n^k \le \lambda_{\max}({\bf V}_n)$.$\Box$
	\vspace{0.5cm}

	With the notation $B_r=\{\boldsymbol{\beta} \ | \  \|\boldsymbol{\beta} - \boldsymbol{\beta}_0\| \le r \}$, $r>0$, we now state a (general) strong consistency theorem.
	\begin{theorem}
		\label{strong-consistency-th}
		Assume that the following conditions hold
		\begin{eqnarray*}
			& (I) & \lambda_{\rm min}[{\bf V}_n] \to \infty, \mbox{ a.s.} \\
			& (S) & \mbox{There exist a constant} \ \boldsymbol{\delta}>0, \mbox{ a random variable } c_0>0 \mbox{ a.s., }\\
			&& \mbox{a random integer } n_1 \geq 1, \mbox{ and a random number } r_1>0 \mbox{ a.s., } \ \mbox{such that} \\
			& & (i) \ |\lambda^T \cD_{n}(\boldsymbol{\beta}) \lambda| >0, \ \mbox{for all} \
			\lambda, \|\lambda\|=1, \ \mbox{and for all} \ \boldsymbol{\beta} \in B_{r_1}, n \geq n_1, \mbox{ a.s.};  \\
			& & (ii) \ \ \lim_{r \to 0} \limsup_{n \to \infty}
			[\lambda_{\max}({\bf V}_n)]^{-1/2-\boldsymbol{\delta}}\sup_{\boldsymbol{\beta} \in
				B_r}\|| \cD_{n}(\boldsymbol{\beta})-\cD_{n}  \||=0 \mbox{ a.s.}; \\
			& & (iii) \ [\lambda_{\max}({\bf V}_n)]^{-1/2-\boldsymbol{\delta}}|\lambda^T\cD_n \lambda| \geq c_0,
			\ \mbox{for all} \ \lambda, \|\lambda\|=1, \ \mbox{and} \ n
			\geq n_1, \mbox{ a.s.}
		\end{eqnarray*}
		\noindent Then, there exists a sequence of random variables $\{ \hat{\boldsymbol{\beta}}_n \}_n$ with values in $\cT$ and a random number $n_0$ such that
		
		(a) $P({\bf q}_n(\hat{\boldsymbol{\beta}}_n)=0, \ \mbox{for all} \ n \geq n_0)=1$;
		
		(b) $\hat{\boldsymbol{\beta}}_n \to \boldsymbol{\beta}_0$ a.s.
	\end{theorem}

	{\bf Proof.} The proof is similar to that of Theorem 4.2 in Balan et al. (2010). We note that here, the normalizer $\alpha_n=\lambda_{\max}({\bf V}_n)$ is random and so Lemma \ref{random-norm-martingale-convergence} should be applied to prove 
	\begin{equation}
	[\lambda_{\max}({\bf V}_n)]^{-1/2-\boldsymbol{\delta}}\|{\bf q}_n\| \to 0, \mbox{ a.s.} \Box \nonumber
	\end{equation}
	
	We now examine the strong consistency of our GEE.
	
	Let ${\bf H}_n'=\sum_{i=1}^n {\bf X}_i^T{\bf A}_i{\bf X}_i$ and note that in the case of deterministic regressors, ${\bf H}_n'={\bf H}_n^{ind}$.
	
	\indent We introduce the following assumptions
	\begin{eqnarray*}
		(D)&& \ \lambda_{\min}({\bf H}_n') \to \infty, \mbox{ a.s. }\\
		(E') && \ \textrm{there exist constants} \ 0<K \le C \ \textrm{such that} \\
		&& \  K \le \lambda_{\min}(\mathcal{R}_{n}^*) \le  \lambda_{\max}(\mathcal{R}_{n}^*) \leq C, \ \textrm{a.s.} \ \forall n \geq 1, \\
		(S_{\boldsymbol{\delta}}')&& \ \textrm{there exist a random integer } N \ge 1,  \textrm{ a constant } 0<\boldsymbol{\delta} \le \frac{1}{2} \mbox { and } \\
		&&\mbox{ a random number } c_0>0 \textrm{ such that} \\
		&&\lambda_{\min}({\bf H}_n') \geq c_0 [\lambda_{\max}({\bf H}_n')]^{1/2+\boldsymbol{\delta}}, \mbox { a.s. } \forall n \geq N.
	\end{eqnarray*}
	\begin{remark}
		\label{remark1}
		{\rm Assumption $(D)$ is a necessary condition (see (1.5) of Lai and Wei, 1982 and the discussion thereafter). Condition $(E')$  mimics condition $(H')$ imposed on $\overline{{\bf R}}_i^{(c)}$, i.e., the approximating sequence $\mathcal{R}_{i}^*$ is also asymptotically positive definite and bounded. In the presence of $(H')$ and $(E')$, assumptions $(I)$ and $(D)$ are equivalent. Assumption $(S'_{\boldsymbol{\delta}})$ was first introduced in Balan et al. (2010), and is a weaker form of assumption $(S_{\boldsymbol{\delta}})$ of Fahrmeir and Kaufmann, (1985). Its use in replacing assumption $(S)(iii)$ above reflects the fact that the leading term of the derivative comes from the generalized linear model (see Fahrmeir and Kaufmann, 1985).}
	\end{remark}

	\vspace{3mm}
	
	If $\mu$ is three times continuously differentiable, for 
	any $r>0$ and $n \geq 1$, let
	\begin{eqnarray*}
		k_n^{[2]}(r)&=&\sup_{\boldsymbol{\beta} \in B_r}\max_{i \leq n, j \leq
			m_i}\left|\frac{\mu''({\bf x}_{ij}^T \boldsymbol{\beta})}{\mu'({\bf x}_{ij}^T
			\boldsymbol{\beta})} \right|, \quad k_n^{[3]}(r)=\sup_{\boldsymbol{\beta} \in B_r}\max_{i \leq
			n, j \leq m_i}\left|\frac{\mu'''({\bf
				x}_{ij}^T \boldsymbol{\beta})}{\mu'({\bf x}_{ij}^T \boldsymbol{\beta})} \right|.
	\end{eqnarray*}
	As in Xie and Yang, 2003 and Balan et al. (2010), condition $(AH)$ is
	\begin{eqnarray*}
		(AH) & & \mbox{there exist random variables} \ C>0,r_0>0 \mbox{ a.s., such that, a.s., } \\
		& & k_{n}^{[l]}(r) \leq C \ \mbox{for all} \ r \le r_0, \  n \geq 1, \ l=2,3.
	\end{eqnarray*}
	For all $n \ge 1$, $r>0$, we define the following random variables
	\begin{eqnarray*}
		\gamma_n^{ '}& =& \max_{i \leq n, j \leq m_i} {\bf
			x}_{ij}^T ({\bf H}_n')^{-1} {\bf x}_{ij}, \
		a_n' = \lambda_{\rm max}({\bf H}_n') \gamma_n^{ '}, \\
		\eta_n(r) &=& \sup_{\boldsymbol{\beta} \in B_r} \max_{i
			\leq n, j\le m_i} \left |\left[\frac{\mu'({\bf x}_{ij}^T \boldsymbol{\beta}')}{\mu'({\bf x}_{ij}^T \boldsymbol{\beta})} \right]^{1/2} - 1 \right|, \\
		\pi_n(r) &=& \sup_{\boldsymbol{\beta} \in B_r} \max_{i
			\leq n} \lambda_{\max}
		[(\cR_{i-1}^*)^{1/2}\cR_{i-1}^*(\boldsymbol{\beta})^{-1}(\cR_{i-1}^*)^{1/2}],
		\\
		d_n(r)&=& \sup_{\boldsymbol{\beta} \in B_r} \max_{i \leq n,l \leq p} \max_{j \le m_i}  \left| \lambda_j \left(\frac{\partial\cR_{i-1}^*(\boldsymbol{\beta})}{\partial \boldsymbol{\beta}_l}
		\right)\right| .
	\end{eqnarray*}
	For $\displaystyle{0<\boldsymbol{\delta} \le \frac{1}{2}}$, we introduce below a set of five conditions, which we label $(C)$. We note that $(C_3')-(C_5)$ compare the a.s. asymptotic behaviours of $\cR_n^*(\boldsymbol{\beta})$ and ${\bf H}_n'$, as $\boldsymbol{\beta}$ approaches the true parameter and $n$ converges to infinity.           
	\begin{eqnarray*}
		(\gamma H') & & \limsup_{n \to \infty}(\gamma_n')^{1/2}[\lambda_{\max}({\bf H}_n')]^{1-\boldsymbol{\delta}} < \infty  \ a.s., \\
		(\pi) & & \lim_{r \to 0} \limsup_{n \to \infty} \pi_n(r) < \infty \  a.s., \\
		(C_3') & & \lim_{r \to 0} \lim_{n \to \infty} r  d_n(r) [\lambda_{\rm max}({\bf H}_n')]^{1/2-\boldsymbol{\delta}}=0 \ a.s., \\
		(C_4) & & \limsup_{n \to \infty }n E[\pi_n^2(r) \tilde a_n'
		\lambda_{\rm max}({\bf
			H}_n')] <\infty,  \ where \ \tilde a_n'=\max\{a_n',
		(a_n')^{2}\},  \\
		(C_5) & & \limsup_{n \to \infty} n E[\pi_n^4(r)d_n^2(r)
		\lambda_{\rm max}({\bf H}_n')]<\infty, \ \mbox{for some} \ r>0.
	\end{eqnarray*} 
	
	\begin{remark}
		\label{NEWconditions}
		{\rm Assume that $(\pi)$ and $(E')$ hold. Then $(C_1)$ and $(C_3)$ of Balan et al. (2010) simplify to $(\gamma H')$ and $(C_3')$, respectively. Furthermore, by Taylor's expansion $(C_3')$ implies $(C_2)$ of Balan et al. (2010), which can thus be replaced by the weaker condition $(\pi)$.}
	\end{remark}
	
	The following result is not stated in Balan et al. (2010).
	\begin{prop}
		\label{conditions_derivative}
		If $(AH)$, $(E')$ and $(C)$ hold, then $(S)(ii)$ holds, i.e.
		\begin{equation}
		\lim_{r \to 0} \limsup_{n \to \infty}
		[\lambda_{\max}({\bf H}_n')]^{-1/2-\boldsymbol{\delta}}\sup_{\boldsymbol{\beta} \in
			B_r}\|| \cD_{n}(\boldsymbol{\beta})-\cD_{n}  \||=0, \mbox{ a.s. } \nonumber
		\end{equation}
	\end{prop}
	{\bf Proof.} The conclusion follows from the proofs of Lemmas 4.6 - 4.9 in Balan et al. (2010), under the new conditions, as shown in Remark \ref{NEWconditions}. Note that the normalizing factor is random, that condition $(K)$ in Balan et al. (2010) can be dropped and that conditions $(E)$ and $(R')$ in Balan et al. (2010) are covered here by condition $(E')$.
	
	\begin{theorem}
		\label{consistency} Assume that conditions $(AH)$, $(H')$, $(D)$, $(E')$ and $(S'_{\boldsymbol{\delta}})$ hold, along with $(C)$.\\
		Then there exists a sequence $\{\widehat{\boldsymbol{\beta}}_n\}_n \subset \mathcal{T}$ and a random number $n_0$ such that 
		
		(a) $P({\bf g}^*_n(\widehat{\boldsymbol{\beta}}_n)=0, \ \textrm{for all} \ n \geq n_0)=1$;
		
		(b) $\widehat{\boldsymbol{\beta}}_n \to \boldsymbol{\beta}_0$ a.s.
	\end{theorem}
	{\bf Proof.} We verify the conditions of Theorem \ref{strong-consistency-th}. By Remark \ref{remark1}, $(I)$ holds. By Proposition \ref{conditions_derivative}, $S(ii)$ holds, $S(i)$ and $S(iii)$ hold as in the proof of Theorem 4.11 in Balan et al. (2010), page 109. $\Box$
	
	\vspace{3mm}
	
	An examination of the proof of Theorem \ref{consistency} gives the following result for the longitudinal generalized linear model, which corresponds to the estimating function ${\bf g}_n^{\rm indep}(\boldsymbol {\beta}),$ with $\boldsymbol{\beta} \in \cT.$
	
	\begin{example}
		Assume that $(H'),$ $(D)$ and $(S_{\delta}')$ hold. Assume further that the following condition holds
		\begin{equation}
		(K') \ \lim_{r \to 0} \limsup_{n \to \infty} \eta_n(r) =0. \nonumber
		\end{equation}
		There exists then a sequence $\{\widehat{\boldsymbol{\beta}}_n\}_{n \ge 1} \subset \cT$ and a random number $n_0$ such that
		
		(a) $P({\bf g}^{\rm indep}_n(\widehat{\boldsymbol{\beta}}_n)=0, \ \textrm{for all} \ n \geq n_0)=1$
		
		(b) $\widehat{\boldsymbol{\beta}}_n \to \boldsymbol{\beta}_0$ a.s.
		
	\end{example}

	{\center \bf Acknowledgements.}
	The authors were funded by a research grant of I. Schiopu-Kratina from the Natural Sciences and Engineering Research Council of Canada.\\

\end{document}